\documentclass{amsart}

\usepackage{amssymb,amsmath,latexsym,amsthm}
\usepackage{pb-diagram}

\newtheorem{theorem}{Theorem}[section]
\newtheorem{lemma}[theorem]{Lemma}

\newtheorem{proposition}[theorem]{Proposition}
\newtheorem{corollary}[theorem]{Corollary}
\newtheorem{conjecture}[theorem]{Conjecture}

\theoremstyle{definition}
\newtheorem*{definition}{Definition}

\newtheorem*{example}{Example}

\newcommand {\Q}{{\mathbb{Q}}}

\newcommand {\Z}{{\mathbb{Z}}}

\newcommand {\F}{{\mathbb{F}}}

\newcommand {\OO}{{\mathcal{O}}}
\newcommand {\D}{{\mathcal{D}}}

\newcommand{\Gal}      {\mathop{\rm {Gal}}}
\newcommand{\Cl}      {\mathop{\rm {Cl}}}
\newcommand{\Norm}        {{\mathcal N}}
\newcommand{\ind}        {{\mathop{\rm ind}}}
\newcommand{\rk}        {{\mathop{\rm rk}}}
\newcommand {\idp}{{\mathfrak{p}}}
\newcommand {\idP}{{\mathfrak{P}}}
\newcommand {\ida}{{\mathfrak{a}}}
\newcommand {\idb}{{\mathfrak{b}}}
\newcommand {\idc}{{\mathfrak{c}}}

\newcommand{\eps}        {\epsilon}
\begin{document}

\title[Asymptotics of number fields]{Asymptotics of number fields and the Cohen--Lenstra heuristics}


\author[J\"urgen {\sc Kl\"uners}]{{\sc J\"urgen} Kl\"uners}
\address{J\"urgen {\sc Kl\"uners}\\
Universit\"at Kassel,\\
Fachbereich f\"ur Mathematik und Informatik\\
Heinrich-Plett-Str. 40\\
34132 Kassel, Germany}
\email{klueners@mathematik.uni-kassel.de}
\begin{abstract}
  We study the asymptotics conjecture of Malle for dihedral groups
  $D_\ell$ of order $2\ell$, where $\ell$ is an odd prime.  We prove
  the expected lower bound for those groups.  For the upper bounds we
  show that there is a connection to class groups of quadratic number
  fields. The asymptotic behavior of those class groups is predicted
  by the Cohen--Lenstra heuristics. Under the assumption of this
  heuristic we are able to prove the expected upper bounds.
\end{abstract}

\maketitle

Dedicated to Michael Pohst on the occasion of his 60th birthday


\bigskip
\section{Introduction}
Let $k$ be a number field and $G\leq S_n$ be a transitive permutation group on $n$ letters. We say
that a finite extension $K/k$ has Galois group $G$ if the normal closure $\hat{K}$ of $K/k$ has
Galois group isomorphic to $G$ and $K$ is the fixed field in $\hat{K}$ under a point stabilizer.
By abuse of notation we write $\Gal(K/k)=G$ in this situation. We let
$$Z(k,G;x):=\#\left\{K/k : \Gal(K/k)=G,\ \Norm_{k/\Q}(d_{K/k})\le
  x\right\}$$
be the number of field extensions of $k$ (inside a fixed
algebraic closure $\bar\Q$) of relative degree~$n$ with Galois group
permutation isomorphic to $G$ (as explained above) and norm of the
discriminant $d_{K/k}$ bounded above by $x$. It is well known that the
number of extensions of $k$ with bounded norm of the discriminant is
finite, hence $Z(k,G;x)$ is finite for all $G$, $k$ and $x\ge1$. We are interested in the
asymptotic behavior of this function for $x\rightarrow\infty$. Gunter Malle \cite{Ma4,Ma5}
has given a precise conjecture how this asymptotic should look like. Before we can state it
we need to introduce some group theoretic definitions.
\begin{definition}
  Let $1\ne G\leq S_n$ be a transitive subgroup acting on $\Omega=\{1,\ldots,n\}$
  and $g\in G$. Then
  \begin{enumerate}
  \item The index $\ind(g):= n- \mbox{ the number of orbits of $g$ on }\Omega.$
  \item $\ind(G):=\min\{\ind(g): 1\ne g\in G\}.$
  \item $a(G):=\ind(G)^{-1}$.
  \item Let $C$ be a conjugacy class of $G$ and $g\in C$. Then $\ind(C):=\ind(g)$.
  \end{enumerate}
\end{definition} 
The last definition is independent of the choice of $g$ since all
elements in a conjugacy class have the same cycle shape. We define an
action of the absolute Galois group of $k$ on the
$\bar{\Q}$-characters of $G$.  The orbits under this action are called
$k$--conjugacy classes. We remark that we get the ordinary conjugacy
classes when $k$ contains all $m$-th roots of unity for $m=|G|$.
\begin{definition}
  For a number field $k$ and a transitive subgroup $1\ne G\leq S_n$ we define:
  $$b(k,G):=\#\{C : C\; k\mbox{-conjugacy class of minimal index }\ind(G)\}.$$
\end{definition}
Now we can state the conjecture of Malle \cite{Ma5}, where we write $f(x) \sim g(x)$ for 
$\lim_{x\rightarrow\infty} \frac{f(x)}{g(x)} =1$.

\begin{conjecture}\label{con}(Malle)
  For all number fields $k$ and all transitive permutation groups $1\ne G$ there exists a 
  constant $c(k,G)>0$ such that
  $$Z(k,G;x) \sim c(k,G)x^{a(G)} \log(x)^{b(k,G)-1},$$
  where $a(G)$ and $b(k,G)$ are given as above.  
\end{conjecture}
We remark that at the time when the conjecture was stated it was only
known to hold for all Abelian groups and the groups $S_3$ and $D_4\leq S_4$.
\begin{example}
  Let $\ell$ be an odd prime and $D_\ell\leq S_\ell$ be the dihedral group of
  order $2\ell$. In this case the non-trivial elements of $D_\ell$ are of order
  2 or $\ell$. In the latter case the index is $\ell-1$. Elements of order 2 
  have 1 fixed point and therefore the index is $(\ell-1)/2$. All elements of
  order 2 are conjugated. Therefore we get:
  $$a(D_\ell)=\frac{2}{\ell-1}\mbox{ and }b(k,D_\ell)=1\mbox{ for all }k.$$
  With the same arguments we get for $D_\ell(2\ell)\leq S_{2\ell}$:
  $$a(D_\ell(2\ell))=\frac{1}{\ell}\mbox{ and }b(k,D_\ell(2\ell))=1\mbox{ for all }k.$$
\end{example}
In \cite{Kl5} we have given a
counter examples to the conjecture. In these counter examples
the $\log$-factor is bigger than expected when certain
subfields of cyclotomic extensions occur as intermediate fields.
Nevertheless the main philosophy of this conjecture should still be true.

The goal of this paper is to prove the conjectured lower bounds for dihedral
groups $D_\ell$, where $\ell$ is an odd prime. We are able to prove the conjectured
upper bounds when we assume a weak version of the Cohen--Lenstra heuristics \cite{CoLe}.
Furthermore we show that the Cohen--Lenstra heuristics is wrong when the conjectured
upper bounds for dihedral groups is broken.

\section{Dihedral groups of order $2\ell$}
\label{sec:dih}

In this section we collect some results about dihedral groups of order
$2\ell$, where $\ell$ is an odd prime. Let $k$ be a number field,
$M/k$ be an extension of degree 2, $N/M$
be an extension of degree $\ell$ such that $N/k$ is normal with Galois group
$D_\ell$. Denote by $K/k$ a subfield of degree $\ell$ and define $d:=\frac{\ell-1}{2}$.
Then we get the following discriminant relation \cite{FiKl}:
\begin{equation}
  \label{eq:disc}
  d_{K/k} = d_{M/k}^d \Norm_{M/k}(d_{N/M})^{1/2}.
\end{equation}
$$\begin{diagram}
\node[3]{N} \arrow[1]{wsw,l,-}{2} \arrow[3]{s,l,-}{D_\ell(2\ell)} \arrow{sse,l,-}{C_\ell} \\
\node[1]{K} \arrow[2]{se,r,-}{D_\ell} \\
\node[4]{M} \arrow{sw,r,-}{2} \\
\node[3]{k} 
\end{diagram}$$
Using this discriminant relation the same approach will work for the group 
$D_\ell(2\ell)\leq S_{2\ell}$ and $D_\ell\leq S_\ell$. The reason is that
up to isomorphy there exists a unique field $K/k$ with Galois group $D_\ell$ 
corresponding to (contained in) a given $N/k$ with Galois group $D_{\ell}(2\ell)$.

Denote by $\D_{k,\ell}$ the set of $D_\ell(2\ell)$--extensions
of $k$ and by $I_k$ the set of ideals of $k$. Then we can define the following mapping:
$$\Psi: \D_{k,\ell} \rightarrow I_k^2, N/k \mapsto \left(d_{M/k}, 
\Norm_{M/k}(d_{N/M})^{1/(2(\ell-1))}\right).$$
In order to define the $2(\ell-1)$th root properly we need the following lemma:
\begin{lemma}
  $\Norm_{M/k}(d_{N/M})$ is a $2(\ell-1)$th power in $I_k$.
\end{lemma}
\begin{proof}
  $N/M$ is a cyclic extension of prime degree. Therefore we have the 
  relation $d_{N/M} = \idc_{N/M}^{\ell-1}$, where $\idc_{N/M}$ is the conductor
  of $N/M$. Let $\idp\subseteq \OO_k$ be a prime ideal dividing $\Norm_{M/k}(d_{N/M})$.
  When $\idp$ is split in $M$ then both prime ideals must divide $\idc_{N/M}$ which
  shows that $\idp^{2(\ell-1)} \mid \Norm_{M/k}(d_{N/M})$. When $\idp$ is inert,
  i.e. $\idp\OO_M=\idP$, we have $\Norm_{M/k}(\idP)=\idp^2$ and we get the desired
  result. The last case is when $\idp\OO_M = \idP^2$. In this case $\idP$ is wildly ramified
  and therefore $\idP^2 \mid \idc_{N/M}$.
  \end{proof}

Let $N\in \D_{k,\ell}$ and $K$ be one of the subfields of $N$ of degree $\ell$. 
Using $d_{N/k}=d_{M/k}^\ell\Norm_{M/k}(d_{N/M})$ and equation \eqref{eq:disc} we get:
\begin{equation}
  \label{eq:disc2}
d_{N/k}=\ida^\ell\idb^{2(\ell-1)} \mbox{ and }d_{K/k}=\ida^d \idb^\ell,  
\end{equation}
where $\Psi(N) =(\ida,\idb)$.

Clearly not every element of $I_k^2$ is in the image of $\Psi$. Let
$(\ida,\idb)\in I_k^2$ be in the image of $\Psi$. Then $\ida$ is
squarefree when we ignore prime ideals lying over $2$.  In other words
$\idp^2\mid \ida$ implies that the prime ideal $\idp$ contains 2. The
ideal $\idb$ is squarefree when we ignore prime ideals lying above
$\ell$. Furthermore, the greatest common divisor of $\ida$ and $\idb$
is only divisible by prime ideals lying above $\ell$. We remark that
these statements are easy to prove, e.g. see \cite[Prop.
10.1.26]{Coh2}. Furthermore in the case $k=\Q$ it can be proven
that $\idb$ itself must be squarefree.

The idea of our proof will be to count elements of $I_k^2$ with the above properties. Unfortunately
$\Psi$ is not injective. In the following we give upper estimates for the number of quadratic
extensions $M/k$ corresponding to a given ideal $\ida$. In a second step we will give
upper estimates for the number of extensions $N/M$ which correspond to the pair $(\ida,\idb)$ when
we assume that $M/k$ is fixed. Altogether we get an upper estimate for the number of elements
in $\D_{k,\ell}$ which map to the same pair $(\ida,\idb)$.

The following lemma is an easy consequence of \cite[Theorem 5.2.9]{Coh2}.
\begin{lemma}
  Let $k$ be a number field and $\ida\subseteq\OO_k$ be an ideal. Denote by $r_u$ the
  unit rank of $k$ and by $r_c$ the 2--rank of the class group $\Cl_k$.
  Then there are at most
  $2^{r_c+r_u+1}$ quadratic extensions $M/k$ which have discriminant $d_{M/k}=\ida$.
\end{lemma}
We remark that in the case $k=\Q$ we get the upper estimate 2 for the number of extensions with 
the same discriminant. We remark that this is the right answer for ideals of the form 
$8\cdot a\cdot \Z$, where $a$ is odd and squarefree since in this case we have the extensions
$\Q(\sqrt{2a})$ and $\Q(\sqrt{-2a})$ which correspond to our ideal.

\begin{lemma}\label{lem:num}
  Let $M/k$ be a quadratic extension such that the $\ell$--rank of the
  class group of $M$ is $r$. Then the number of
  $D_\ell(2\ell)$--extensions $N/M/k$ which are subfields of the ray
  class field of $\idb\OO_M$ is bounded from above by
  $\frac{\ell^{r+\omega(\idb)}-1}{\ell-1}$, where $\omega(\idb)$ denotes the
  number of different prime factors of $\idb$.
\end{lemma}
\begin{proof}
  Define $H:=\Gal(M/k)=\langle \sigma \rangle$ and consider the
  ray class group of $\idb\OO_M$ modulo $\ell$th powers as an $H$-module.
  In \cite[Sections 5 and 6]{FiKl} it is proved that
  $D_\ell$-extensions correspond to invariant subspaces of dimension 1
  and eigenvalue $-1$.  Let $\idp$ be a prime ideal of $\OO_k$ not
  lying over $\ell$. We distinguish two cases. First we assume that
  $\idp$ decomposes in $M$. In this case the $\ell$-rank of the
  residue class ring $\OO_M/\idb\OO_M$ is 2 or 0 depending on $\ell
  \mid |(\OO_k/\idp)^\times|$ or not.  If the rank is 2 then $\sigma$
  interchanges the two prime ideals lying over $\idp$ and we find a
  1-dimensional subspace with eigenvalue 1 and one with eigenvalue -1.
  In case that $\idp$ is inert the $\ell$-rank of the ray class
  group increases by at most 1. Similar arguments show the same
  result for prime ideals $\idp$ lying over $\ell$.
\end{proof}
We remark that in the above lemma we give an upper estimate for all
$C_\ell$--extensions $N/M$ such that the conductor is a divisor of
$\idb$.

Altogether we get the following upper bound for the number of fibers
of $(\ida,\idb)$ under $\Psi$:
\begin{equation}
  \label{eq:upp}
2^{\rk_2(\Cl_k)+r_u+1} \cdot
\frac{\ell^{r+\omega(\idb)}-1}{\ell-1},  
\end{equation}
where $r$ is the maximal
$\ell$--rank of the class group of a quadratic extension $M/k$ of
discriminant $\ida$ and $r_u$ is the unit rank of $k$. Since our ground field $k$ is fixed, we have that
$2^{\rk_2(\Cl_k)+r_u+1}$ is a constant depending on $k$. We will see
that $\omega(\idb)$ will cause no problems. The critical part in this
estimate is the dependency on the class group of a field $M$ which we
do not explicitly know.

In order to get good upper estimates we need to control the $\ell$--rank
of quadratic extensions. In order to simplify the situation let us restrict
to the case $k=\Q$. In the general case we get similar results when we
assume the corresponding things for relative quadratic extensions $M/k$.
The Cohen--Lenstra--heuristic predicts the behavior of the class group
of quadratic number fields. Let us state a special case of this conjecture
\cite[C6 and C10]{CoLe} which is only proven for $\ell=3$. In the following
conjecture all sums are over fundamental discriminants $D$. We define
$r_D:=\rk_\ell(\Cl_{\Q(\sqrt{D})})$.

\begin{conjecture} (Cohen-Lenstra) \label{coh-len}

  The average of $\ell^{r_D}-1$ over all imaginary quadratic fields is 1, 
  i.e.
  $$\frac{\sum_{-D\leq x} \ell^{r_D}-1}
{\sum_{-D\leq x} 1} \longrightarrow 1 \mbox{ for }x\longrightarrow \infty,$$
where only fundamental discriminants are considered in the sums.
  The average of $\ell^{r_D}-1$ over all real quadratic fields is $\ell^{-1}$.
\end{conjecture}

It is well known that $Z(\Q,C_2;x)\sim c(C_2) x$, where $c(C_2)$ is explicitly
known, see e.g. \cite{CoDiOl2}. For our purposes it will be enough to
assume the following, where $O$ is the usual Landau symbol.
\begin{equation}
  \label{eq:CoLe_i}
  \sum_{-D\leq x} \ell^{r_D} = O(x),
\end{equation}
\begin{equation}
  \label{eq:CoLe_r}
  \sum_{D\leq x} \ell^{r_D} = O(x).
\end{equation}

\begin{theorem}\label{main:upper}
  Assume equations \eqref{eq:CoLe_i} and \eqref{eq:CoLe_r}. Then for all
  odd primes $\ell$ we find constants $c(\ell)$ and $\hat c(\ell)$ such that:
$$Z(\Q,D_\ell;x) \leq c(\ell) x^{1/d} =c(\ell) x^{a(D_\ell)}, 
\mbox{ where }d=\frac{\ell-1}{2},$$
$$Z(\Q,D_\ell(2\ell);x) \leq \hat{c}(\ell) x^{1/\ell}=\hat{c}(\ell) x^{a(D_\ell(2\ell))}.$$
\end{theorem}
\begin{proof}
  Since $\Z$ is a principal ideal domain we can parameterize ideals in $\Z$ by positive
  integers. Using equations \eqref{eq:disc2} and \eqref{eq:upp} we get:
  $$Z(\Q,D_\ell(2\ell),x)\leq \sum_{D^\ell b^{2(\ell-1)}\leq x} 
  \frac{\ell^{\omega(b)+r_D}-1}{\ell-1} +
\sum_{(-D)^\ell b^{2(\ell-1)}\leq x}  \frac{\ell^{\omega(b)+r_D}-1}{\ell-1},
$$
where $D$ runs over the fundamental discriminants, positive in the first sum and negative
in the second one. Let us restrict to the first sum. The other one will be estimated in the
same way.
$$
\sum_{D^\ell b^{2(\ell-1)}\leq x} 
\frac{\ell^{\omega(b)+r_D}-1}{\ell-1} \leq
\sum_{D^\ell b^{2(\ell-1)}\leq x} \ell^{\omega(b)} \ell^{r_D}$$
$$= \sum_{b^{2(\ell-1)}\leq x} \ell^{\omega(b)} 
\sum_{D\leq \frac{x^{1/\ell}}{b^{2(\ell-1)/\ell}}} \ell^{r_D}
\stackrel{\eqref{eq:CoLe_r}}{\leq} 
\tilde{c}\sum_{b^{2(\ell-1)}\leq x} \frac{\ell^{\omega(b)}x^{1/\ell}}{b^{2(\ell-1)/\ell}}
$$
$$=\tilde{c}x^{1/\ell} \sum_{b^{2(\ell-1)}\leq x} \frac{\ell^{\omega(b)}}{b^{2(\ell-1)/\ell}}
=c x^{1/\ell}.$$ The last sum converges since $2(\ell-1)/\ell>1$.
The result for $Z(\Q,D_\ell,x)$ follows in the same way when we sum over 
$D^d b^{\ell-1}\leq x$.
\end{proof}

As already remarked we can prove similar results for arbitrary ground fields $k$ when
we assume corresponding results for the $\ell$--rank of the class groups of quadratic
extensions $M/k$. This becomes quite technical which is the reason that we do not give
this case here.

Conjecture \ref{coh-len} is proved for $\ell=3$ (\cite{DaWr}) 
which gives the following corollary.
\begin{corollary}
  \begin{enumerate}
  \item   $Z(\Q,D_3;x) \leq c(\ell) x^{a(D_3)}$.
  \item   $Z(\Q,D_3(6);x) \leq c(\ell) x^{a(D_3(6))}$.
  \end{enumerate}
\end{corollary}
We remark that in the paper \cite{DaWr} the stronger result:
$Z(k,S_3;x) \sim c(k,S_3) x$ is obtained.

If we want to have unconditional upper bounds for dihedral groups, the best thing
we can do at the moment is the following:
$$\ell^r_D \leq \#\Cl\nolimits_{\Q(\sqrt{D})} =O(D^{1/2}\log(D)).$$
If we use this estimate and use the method of the proof of Theorem \ref{main:upper}
we can prove the following.

\begin{theorem}\label{main:upper2}
  Let $\ell$ be an odd prime. Then for all $\eps>0$ we can 
  find constants $c(\ell,\eps)$ and $\hat c(\ell,\eps)$ such that:
$$Z(\Q,D_\ell;x) \leq c(\ell,\eps) x^{3a(D_\ell)/2+\eps}, $$
$$Z(\Q,D_\ell(2\ell);x) \leq \hat{c}(\ell,\eps) x^{3a(D_\ell(2\ell))/2+\eps}.$$
\end{theorem}

\section{Lower bounds}
\label{sec:lower}

In this section we are interested in lower bounds for the number of fields with
dihedral Galois group. First we show that Cohen--Lenstra heuristic also
provides lower bounds for our asymptotics. We prove that the non-validity
of equations \eqref{eq:CoLe_i} or \eqref{eq:CoLe_r} implies that the asymptotics
conjecture for the corresponding dihedral groups is wrong as well. This shows
that the class group of intermediate fields is an important obstruction. As in the 
preceding section we assume $k=\Q$ to simplify everything. The following lemma
is well known.
\begin{lemma}\label{unram}
  Let $M/\Q$ be a quadratic extension and $N/M$ be a cyclic unramified extension of degree
  $\ell$, where $\ell>2$ is prime. Then $\Gal(N/\Q)=D_\ell(2\ell)$.
\end{lemma}
\begin{proof}
  $\Gal(N/\Q)$ is a transitive subgroup of the wreath product $C_\ell \wr C_2$ and
  therefore one of the following groups: $C_\ell,D_\ell(2\ell),$ or $C_\ell \wr C_2$.
  Furthermore $\Gal(N/\Q)$ is generated by its inertia groups. All ramified primes
  have order 2 and the dihedral group is the unique group in this list which is
  generated by elements of order 2.
\end{proof}

Now we count dihedral extensions which are unramified over their quadratic subfield.
$$Y(\Q,D_\ell(2\ell);x):=$$
$$\#\{N/\Q:\Gal(N/\Q)=D_\ell(2\ell),\Norm(d_{N/\Q})\leq x,
N/M \mbox{ unramified}\}.$$
When we define $Y(\Q,D_\ell;x)$ in an analogous way we get:
\begin{theorem}
  Assume Conjecture \ref{coh-len}. Then 
  \begin{enumerate}
  \item $Y(\Q,D_\ell(2\ell);x) \sim c x^{1/\ell}= c x^{a(D_\ell(2\ell))}.$
  \item $Y(\Q,D_\ell;x) \sim \tilde c x^{2/(\ell-1)} = \tilde{c} x^{a(D_\ell)}.$
  \end{enumerate}
\end{theorem}
\begin{proof}
  As in the proof of Theorem \ref{main:upper} we consider real and complex quadratic
  fields separately. 
  $$Y(\Q,D_\ell(2\ell);x) = \sum_{D^\ell\leq x} \frac{\ell^{r_D}-1}{\ell-1} 
  +\sum_{-D^\ell\leq x} \frac{\ell^{r_D}-1}{\ell-1}$$
$$=\frac{1}{\ell-1} \left(\sum_{D\leq x^{1/\ell}} \ell^{r_D}-1 + 
    \sum_{-D\leq x^{1/\ell}} \ell^{r_D}-1\right).$$
  Using conjecture \ref{coh-len} we get the formula in our theorem. The second
  one can be proved in an analogous way.
\end{proof}
Now assume that one of the equations \eqref{eq:CoLe_i} or \eqref{eq:CoLe_r} is wrong.
Using the same arguments as in the last proof we get a lower bound for dihedral
groups $D_\ell$ which contradicts the asymptotics conjecture \ref{con}.

In the last part of this paper we prove the lower bound for dihedral groups
unconditionally. We need the following result which is a special case of 
\cite[Theorem 4.2]{DaWr}.
\begin{proposition}\label{prop}
  Let $\ell$ be an odd prime and $p,q\equiv 1 \bmod \ell$ be two primes. Then the
  number of quadratic extensions $M/\Q$ which are split in $p$ and $q$ grows
  asymptotically like $cx$ for some explicit constant $c$.
\end{proposition}
The quadratic extensions given in the last proposition have the nice property that
they admit a dihedral extension with bounded discriminant.
\begin{lemma}\label{lem:ex}
  Let $M/\Q$ be a quadratic extension which is split in two primes $p,q$ which are
  congruent to $1\bmod \ell$ for an odd prime $\ell$. Then there exists an extension
  $N/M$ which is at most ramified in primes lying above $p$ and $q$ such that
  $\Gal(N/\Q)=D_\ell(2\ell)$.
\end{lemma}
\begin{proof}
  Using Lemma \ref{unram} we can assume that $\ell\nmid \#\Cl_M$. Define $\ida:=pq\OO_M$
  and consider the ray class group $\Cl_\ida$. The multiplicative group of the
  residue ring $\OO_M/\ida$ has $\ell$--rank 4. Furthermore the $\ell$--rank of
  the unit group of $\OO_M$ is bounded above by 1. Using the canonical diagram
  for ray class groups, see e.g. \cite[p. 126]{La}, shows that the $\ell$--rank
  of $\Cl_\ida$ is at least 3. Let $\sigma$ be the generator of the Galois group
  of $M/\Q$ and denote by $A_\ell$ the $\F_\ell[C_2]$--module $\Cl_\ida/\Cl_\ida^\ell$,
  using the canonical action of $\sigma$. Now we can decompose
  $A_\ell = A_\ell^+ \oplus A_\ell^-$, where $A_\ell^+:=\{a \in A_\ell\mid \sigma(a)=a\}$
  and $A_\ell^-:=\{a\in A_\ell\mid \sigma(a)=a^{-1}\}$. 
  Using \cite[Sections 5 and 6]{FiKl} we know that elements of $A_\ell^+$ correspond
  to $C_2\times C_\ell$ extensions and elements of $A_\ell^-$ correspond to
  $D_\ell(2\ell)$ extensions $N/\Q$. Since $C_2\times C_\ell$ extensions of $\Q$ admit
  a subfield with Galois group $C_\ell$ over $\Q$, we see that $\rk_\ell(A_\ell)^+=2$.
  This implies $\rk_\ell(A_\ell^-)\geq 1$ and therefore the existence of our $D_\ell(2\ell)$
  extension.
\end{proof}
Let $N/M$ be such an extension given in the above lemma.
Then 
\begin{equation}
  \label{eq:disc3}
\Norm(d_{N/\Q})\leq \Norm(d_M)^\ell (pq)^{2(\ell-1)} \mbox{ and }
\Norm(d_{K/\Q})\leq \Norm(d_K)^{(\ell-1)/2} (pq)^{\ell-1},  
\end{equation}
  where $K$ is a subfield of $N$ of degree $\ell$.
Now we can prove the following theorem.
\begin{theorem}
  Let $\ell$ be an odd prime. Then there exist positive constants $c_1(\ell),c_2(\ell)$ 
  such that
  \begin{enumerate}
  \item $Z(\Q,D_\ell;x)\geq c_1(\ell) x^{a(D_\ell)} \mbox{ for }x\mbox{ large enough}.$
  \item $Z(\Q,D_\ell(2\ell);x)\geq c_2(\ell) x^{a(D_\ell(2\ell))} \mbox{ for }x\mbox{ large enough}.$
  \end{enumerate}
\end{theorem}
\begin{proof}
  Choose primes $p,q\equiv \bmod \;\ell$. For every quadratic extension $M/\Q$ which
  is split in $p$ and $q$ we find a $D_\ell$ extension by Lemma \ref{lem:ex} which has bounded discriminant
  as given in \eqref{eq:disc3}. Clearly, for different $M$ those $D_\ell$ extensions are
  different.  Using our discriminant formula we get:
  $$Z(\Q,D_\ell;x) \geq \sum_{D^{(\ell-1)/2}(pq)^{\ell-1}\leq x} 1 = 
  \sum_{D\leq (x/(pq)^{\ell-1})^{2/(\ell-1)}} 1,$$
    where we sum only over fundamental discriminants $D$ such that $p$ and $q$ are
    split in $\Q(\sqrt{D})$. Using Proposition \ref{prop} we get that the last sum grows asymptotically
    like $c x^{2/(\ell-1)}$ which proves our first formula. The second one
    can be proved analogously.
\end{proof}

\end{document}